%
\catcode`@=11
%
%
\def\bibn@me{R\'ef\'erences}
\def\bibliographym@rk{\centerline{{\sc\bibn@me}}
	\sectionmark\section{\ignorespaces}{\unskip\bibn@me}
	\bigbreak\bgroup
	\ifx\ninepoint\undefined\relax\else\ninepoint\fi}
%
%
%
\let\refsp@ce=\ 
\let\bibleftm@rk=[
\let\bibrightm@rk=]
%
%
%
\def\numero{n\raise.82ex\hbox{$\fam0\scriptscriptstyle o$}~\ignorespaces}
%
%
\newcount\equationc@unt
\newcount\bibc@unt
\newif\ifref@changes\ref@changesfalse
\newif\ifpageref@changes\ref@changesfalse
\newif\ifbib@changes\bib@changesfalse
\newif\ifref@undefined\ref@undefinedfalse
\newif\ifpageref@undefined\ref@undefinedfalse
\newif\ifbib@undefined\bib@undefinedfalse
\newwrite\@auxout
%
%
\def\eqnum{\global\advance\equationc@unt by 1%
\edef\lastref{\number\equationc@unt}%
\eqno{(\lastref)}}
%
%
%
%
%
%
\def\re@dreferences#1#2{{%
	\re@dreferenceslist{#1}#2,\undefined\@@}}
\def\re@dreferenceslist#1#2,#3\@@{\def\next{#2}%
	\expandafter\ifx\csname#1@@\meaning\next\endcsname\relax
	??\immediate\write16
	{Warning, #1-reference "\next" on page \the\pageno\space
	is undefined.}%
	\global\csname#1@undefinedtrue\endcsname
	\else\csname#1@@\meaning\next\endcsname\fi
	\ifx#3\undefined\relax
	\else,\refsp@ce\re@dreferenceslist{#1}#3\@@\fi}
%
%
%
\def\newlabel#1#2{{\def\next{#1}\newl@bel#2}}
\def\newl@bel#1#2{%
	\expandafter\xdef\csname ref@@\meaning\next\endcsname{#1}%
	\expandafter\xdef\csname pageref@@\meaning\next\endcsname{#2}}
\def\label#1{{%
	\toks0={#1}\message{ref(\lastref) \the\toks0,}%
	\ignorespaces\immediate\write\@auxout%
	{\noexpand\newlabel{\the\toks0}{{\lastref}{\the\pageno}}}%
	\def\next{#1}%
	\expandafter\ifx\csname ref@@\meaning\next\endcsname\lastref%
	\else\global\ref@changestrue\fi%
	\newlabel{#1}{{\lastref}{\the\pageno}}}}
\def\ref#1{\re@dreferences{ref}{#1}}
\def\pageref#1{\re@dreferences{pageref}{#1}}
%
%
\def\bibcite#1#2{{\def\next{#1}%
	\expandafter\xdef\csname bib@@\meaning\next\endcsname{#2}}}
\def\cite#1{\bibleftm@rk\re@dreferences{bib}{#1}\bibrightm@rk}
%
%
\def\beginthebibliography#1{\bibliographym@rk
	\setbox0\hbox{\bibleftm@rk#1\bibrightm@rk\enspace}
	\parindent=\wd0
	\global\bibc@unt=0
	\def\bibitem##1{\global\advance\bibc@unt by 1
		\edef\lastref{\number\bibc@unt}
		{\toks0={##1}
		\message{bib[\lastref] \the\toks0,}%
		\immediate\write\@auxout
		{\noexpand\bibcite{\the\toks0}{\lastref}}}
		\def\next{##1}%
		\expandafter\ifx
		\csname bib@@\meaning\next\endcsname\lastref
		\else\global\bib@changestrue\fi%
		\bibcite{##1}{\lastref}
		\medbreak
		\item{\hfill\bibleftm@rk\lastref\bibrightm@rk}%
		}
	}
\def\endthebibliography{\egroup\par}
%
%
%
\def\@closeaux{\closeout\@auxout
	\ifref@changes\immediate\write16
	{Warning, changes in references.}\fi
	\ifpageref@changes\immediate\write16
	{Warning, changes in page references.}\fi
	\ifbib@changes\immediate\write16
	{Warning, changes in bibliography.}\fi
	\ifref@undefined\immediate\write16
	{Warning, references undefined.}\fi
	\ifpageref@undefined\immediate\write16
	{Warning, page references undefined.}\fi
	\ifbib@undefined\immediate\write16
	{Warning, citations undefined.}\fi}
%
%
\immediate\openin\@auxout=\jobname.aux
\ifeof\@auxout \immediate\write16
  {Creating file \jobname.aux}
\immediate\closein\@auxout
\immediate\openout\@auxout=\jobname.aux
\immediate\write\@auxout {\relax}%
\immediate\closeout\@auxout
\else\immediate\closein\@auxout\fi
%
%
\input\jobname.aux
\immediate\openout\@auxout=\jobname.aux
%
%

\def\bibn@me{R\'ef\'erences bibliographiques}
%
\def\bibliographym@rk{\bgroup}
%
%
\outer\def\bye{ 	\par\vfill\supereject\end}

\def\Q{{\bf {Q}}}

\def\Z{{\bf Z}}

\overfullrule=0pt

\magnification=1200

  \def\pro{\noindent {\bf{Proof : }}}

\def\house#1{\setbox1=\hbox{$\,#1\,$}%
\dimen1=\ht1 \advance\dimen1 by 2pt \dimen2=\dp1 \advance\dimen2 by 2pt
\setbox1=\hbox{\vrule height\dimen1 depth\dimen2\box1\vrule}%
\setbox1=\vbox{\hrule\box1}%
\advance\dimen1 by .4pt \ht1=\dimen1
\advance\dimen2 by .4pt \dp1=\dimen2 \box1\relax}

\def\Card{{\rm Card}}

  \def\eps{{\varepsilon}}

\def\disp{\displaystyle} 
  \def\noi{\noindent}

\def\build#1_#2^#3{\mathrel{\mathop{\kern 0pt#1}\limits_{#2}^{#3}}}

\def\date {le\ {\the\day}\ \ifcase\month\or 
janvier\or fevrier\or mars\or avril\or mai\or juin\or juillet\or
ao\^ut\or septembre\or octobre\or novembre\or 
d\'ecembre\fi\ {\oldstyle\the\year}}

\font\fivegoth=eufm5 \font\sevengoth=eufm7 \font\tengoth=eufm10

\newfam\gothfam \scriptscriptfont\gothfam=\fivegoth
\textfont\gothfam=\tengoth \scriptfont\gothfam=\sevengoth

\def\cqfd{\unskip\kern 6pt\penalty 500 \raise 0pt\hbox{\vrule\vbox 
to6pt{\hrule width 6pt \vfill\hrule}\vrule}\par}

\def\pro{\noindent {\it Proof. }}

\def\smallsquare{\vbox{\hrule\hbox{\vrule height 1 ex\kern 1 ex\vrule}\hrule}}
\def\cqfd{\hfill \smallsquare\vskip 3mm}

\def\Card{{\rm Card}}

\def\uU{{\overline U}}

\def\uW{{\overline W}}

\centerline{}

\vskip 4mm

\centerline{
\bf  Automatic continued fractions are transcendental or quadratic}

\vskip 2mm

\centerline{Les fractions continues automatiques sont transcendantes 
ou quadratiques}

\vskip 13mm

\centerline{
Yann B{\sevenrm UGEAUD}
\footnote{}{\rm
}}

{\narrower\narrower
\vskip 15mm

\proclaim Abstract. {
We establish new combinatorial transcendence criteria
for continued fraction expansions.
Let $\alpha = [0; a_1, a_2, \ldots]$ be an algebraic number 
of degree at least three.
One of our criteria implies that the sequence of partial
quotients $(a_{\ell})_{\ell \ge 1}$ of $\alpha$ 
is not `too simple' (in a suitable sense) and
cannot be generated by a finite automaton.
}

}

{\narrower\narrower
\vskip 10mm

\proclaim R\'esum\'e. {
Nous \'etablissons de nouveaux crit\`eres combinatoires
de transcendance pour des d\'eveloppements en fraction continue.
Soit $\alpha = [0; a_1, a_2, \ldots]$ un nombre alg\'ebrique de
degr\'e au moins \'egal \`a trois. 
L'un de nos crit\`eres entra\^\i ne que la suite $(a_{\ell})_{\ell \ge 1}$
des quotients partiels de $\alpha$ n'est pas trop simple (en un certain sens)
et ne peut pas \^etre engendr\'ee par un automate fini.
}

}

\vskip 11mm

\centerline{\bf 1. Introduction and results}

\vskip 6mm

A well-known open question in Diophantine
approximation asks whether the 
continued fraction expansion of an irrational algebraic number 
$\alpha$ either is ultimately periodic 
(this is the case if, and only if, $\alpha$ is a 
quadratic irrational), or contains arbitrarily large partial quotients. 
As a preliminary step towards its resolution, several 
transcendence criteria for continued fraction expansions have been 
established recently \cite{AdBu05,AdBu07e,AdBu07f,AdBuDa}
(we refer the reader to these papers for references
to earlier works, which include \cite{Mai06,Bak62,Que98,ADQZ}) 
by means of a deep tool
from Diophantine approximation, namely the Schmidt Subspace Theorem
(see Theorem 2.1 below).
In the present note, we show how a slight modification 
of their proofs allows us to considerably improve two of these
criteria. We begin by pointing out two important
consequences of one of our new criteria.
Thus, we solve two problems addressed and discussed in \cite{AdBu05}
and we establish for continued fraction expansions of algebraic numbers
the analogues of the results of \cite{AdBu07d} on the expansions
of algebraic numbers to an integer base.

Throughout this note, ${\cal A}$ denotes a finite or infinite set,
called the alphabet.
We identify a sequence ${\bf a}=(a_\ell)_{\ell \ge 1}$ 
of elements from ${\cal A}$
with the infinite word $a_1 a_2 \ldots a_\ell \ldots$, as well
denoted by ${\bf a}$.
This should not cause any confusion.

For $n \ge 1$, we denote by $p(n, {\bf a})$ the
number of distinct blocks of $n$
consecutive letters occurring in the
word ${\bf a}$, that is,
$$
p(n, {\bf a}) := \Card \{a_{\ell+1} \ldots a_{\ell + n} : \ell \ge 0\}.
$$
The function 
$n \mapsto p(n, {\bf a})$ is called the
{\it complexity function} of ${\bf a}$.
A well-known result of Morse and Hedlund \cite{MoHe38,MoHe40}
asserts that $p(n, {\bf a}) \ge n + 1$ for $n \ge 1$,
unless ${\bf a}$ is ultimately periodic (in which case
there exists a constant $C$ such that 
$p(n, {\bf a}) \le C$ for $n \ge 1$).

\goodbreak

Our first result asserts that the complexity function
of the sequence of partial quotients $(a_{\ell})_{\ell \ge 1}$ of
an algebraic number 
$$
[0; a_1, a_2, \ldots, a_\ell, \ldots] = 
 {1 \over \disp a_1 +
{\strut 1 \over \disp a_2 + {\strut 1 \over \disp \ldots}}}
$$
of degree at least three
cannot increase too slowly.

\proclaim Theorem 1.1.
Let ${\bf a}=(a_\ell)_{\ell \ge 1}$ be a sequence of positive integers
which is not ultimately periodic.
If the real number
$$
[0; a_1, a_2, \ldots,a_\ell, \ldots]
$$ 
is algebraic, then
$$
\lim_{n \to + \infty} \, {p(n, {\bf a}) \over n} = + \infty.   \eqno (1.1)
$$

Theorem 1.1 improves Theorem 7 from \cite{ADQZ}
and Theorem 4 from \cite{AdBu05},
where 
$$
\lim_{n \to + \infty} \, {p(n, {\bf a}) - n} = + \infty
$$
was proved instead of (1.1).
This gives a positive answer to Problem 3 of \cite{AdBu05}
(we have chosen here a different formulation).


An infinite sequence ${\bf a}=(a_\ell)_{\ell \geq 1}$ is 
an automatic sequence if it can be generated by a finite
automaton, that is, if
there exists an integer $k \ge 2$ such that $a_\ell$ is a 
finite-state function of the 
representation of $\ell$ in base $k$, for every $\ell \ge 1$.
We refer the reader to \cite{AlSh} for a more precise
definition and examples of automatic sequences.
Let $b \ge 2$ be an integer.
In 1968, Cobham \cite{Cob68} asked whether a real number
whose $b$-ary expansion can be generated by a finite automaton
is always either rational or transcendental. A positive answer
to Cobham's question
was recently given in \cite{AdBu07d}. 
We addressed in \cite{AdBu05} 
the analogous question for continued fraction expansions.
Since the complexity function of every automatic sequence ${\bf a}$
satisfies $p(n, {\bf a}) = O(n)$ (this was proved
by Cobham \cite{Cob72} in 1972), Theorem 1.1 implies
straightforwardly a negative answer to Problem 1 of \cite{AdBu05}.

\proclaim Theorem 1.2. 
The continued fraction expansion of an algebraic number 
of degree at least three cannot be generated by a
finite automaton.


The proofs of Theorems 1.1 and 1.2 rest ultimately on
a combinatorial transcendence criterion established by
means of the
Schmidt Subspace Theorem.
This is also the case for the
similar results about expansions 
of irrational algebraic numbers to an integer base, see \cite{AdBu07d,AdBuLu04}.

Before stating our criteria, we introduce some notation.
The length of a word
$W$ on the alphabet ${\cal A}$, that is, the number of letters
composing $W$, is denoted by $\vert W\vert$.
We denote the mirror image of a finite word $W := a_1 \ldots a_\ell$ 
by $\uW := a_\ell \ldots a_1$. In particular, $W$ is a palindrome if,
and only if, $W = \uW$.

Let ${\bf a}=(a_\ell)_{\ell \ge 1}$ 
be a sequence of elements from ${\cal A}$.
We say that ${\bf a}$ 
satisfies Condition $(*)$ if ${\bf a}$ is not  
ultimately periodic and if there exist 
three sequences of finite words $(U_n)_{n \ge 1}$,   
$(V_n)_{n \ge 1}$ and $(W_n)_{n \ge 1}$ such that:

\medskip

\item{\rm (i)} For every $n \ge 1$, either the word $W_n U_n V_n U_n$ 
or the word $W_n U_n V_n \uU_n$ is a prefix
of the word ${\bf a}$;

\smallskip

\item{\rm (ii)} The sequence
$({\vert V_n\vert} / {\vert U_n\vert})_{n \ge 1}$ is bounded from above;  

\smallskip

\item{\rm (iii)} The sequence
$({\vert W_n\vert} / {\vert U_n\vert})_{n \ge 1}$ is bounded 
from above;   

\smallskip

\item{\rm (iv)} The sequence $(\vert U_n\vert)_{n \ge 1}$ is 
increasing.

\medskip
\goodbreak

Equivalently, the word ${\bf a}$ satisfies Condition $(*)$ if
there exists a positive real number $\eps$ such that, for arbitrarily
large integers $N$, the prefix $a_1 a_2 \ldots a_N$ of ${\bf a}$
contains two disjoint occurrences of a word of length $[\eps N]$
or it contains a word $U$ of length $[\eps N]$ and its mirror image $\uU$,
provided that $U$ and $\uU$ do not overlap. Here and below, $[ \cdot ]$ denotes the
integer part function.

We summarize our two new 
combinatorial transcendence criteria in the following theorem.

\proclaim Theorem 1.3.
Let ${\bf a}=(a_\ell)_{\ell \ge 1}$ be a sequence of positive integers.
Let $(p_\ell/q_\ell)_{\ell \ge 1}$ denote the sequence of convergents to 
the real number
$$
\alpha:= [0; a_1, a_2, \ldots,a_\ell, \ldots].
$$ 
Assume that the sequence $(q_\ell^{1/\ell})_{\ell \ge 1}$ is bounded. 
If ${\bf a}$ satisfies Condition $(*)$,  
then $\alpha$ is transcendental.

Theorem 1.3 is the disjoint union of two transcendence criteria.
A first one applies to stammering continued fractions, where the terminology 
`stammering' means that in (i) the word $W_n U_n V_n U_n$ is a prefix
of the word ${\bf a}$ for infinitely many $n$; see Theorem 3.1.
A second one is concerned with quasi-palindromic continued fractions,
where the terminology 
`quasi-palindromic' means that in (i) the word $W_n U_n V_n \uU_n$ is a prefix
of the word ${\bf a}$ for infinitely many $n$; see Theorem 5.1.
The condition that the sequence $(q_\ell^{1/\ell})_{\ell \ge 1}$ 
has to be bounded is not very restrictive, since it is
satisfied by almost all real numbers (in the sense of the Lebesgue measure).
Furthermore, it is clearly satisfied when $(a_{\ell})_{\ell \ge 1}$
is bounded. 
Note that this condition can be removed if ${\bf a}$ begins with arbitrarily large
squares $U_n U_n$ (Theorem 2.1 from \cite{AdBuDa})
or with arbitrarily large palindromes $U_n \uU_n$
(Theorem 2.1 from \cite{AdBu07f}).

Theorem 1.3 encompasses all the combinatorial
transcendence criteria for continued
fraction expansions established in \cite{AdBu05,AdBu07e,AdBu07f,AdBuDa}
under the assumption that the sequence $(q_\ell^{1/\ell})_{\ell \ge 1}$ 
is bounded.

Let ${\bf a}$ be a sequence of positive integers.
If there exist three sequences of finite words $(U_n)_{n \ge 1}$,   
$(V_n)_{n \ge 1}$ and $(W_n)_{n \ge 1}$
such that $\limsup_{n \to + \infty} |W_n| / |U_n|$
is sufficiently small and ${\bf a}$ satisfies Condition $(*)$, then the transcendence
of $[0; a_1, a_2, \ldots ]$ was already proved in
\cite{AdBu05,AdBuDa,AdBu07f}. The novelty in Theorem 1.3
is that we allow $|W_n|$ to be large, provided however that the
quotients $|W_n| / |U_n|$ remain bounded independently of $n$.
This is crucial for the proofs of Theorems 1.1 and 1.2.

At present, we do not know any transcendence criterion involving palindromes
for expansions to integer bases; however, see \cite{AdBu06}. 

We end this section with an application of
Theorem 3.1 to quasi-periodic continued fractions.

\proclaim Theorem 1.4.
Consider the quasi-periodic continued fraction 
$$
\alpha =[0 ; a_1,\ldots,a_{n_0-1},
\underbrace{a_{n_0},\ldots,a_{n_0+r_0-1}}_{\lambda_0}
,\underbrace{a_{n_1},\ldots,a_{n_1+r_1-1}}_{\lambda_1},\ldots],
$$
where the notation means that $n_{k+1}=n_k+\lambda_kr_k$ and the 
$\lambda$'s indicate 
the number of times a block of partial quotients is repeated. 
Let $(p_\ell / q_\ell)_{\ell \ge 1}$ denote the sequence of convergents to $\alpha$.
Assume that the sequence $(q_\ell^{1/\ell})_{\ell \ge 1}$ is bounded.
If  the sequence 
$(a_\ell)_{\ell \geq 1}$ is not ultimately periodic and  
$$
\liminf_{k\to\infty} {\lambda_{k+1} \over \lambda_k} > 1, \eqno (1.2)
$$
then the real number $\alpha$ is transcendental.

Theorem 1.4  improves Theorem 3.4 from \cite{AdBu07e}, where, instead of 
the assumption (1.2), the stronger condition
$\liminf_{k\to\infty} {\lambda_{k+1} / \lambda_k} > 2$ was required.

\vskip 7mm

\centerline{\bf 2. Auxiliary results}

\vskip 5mm

We gather below several classical results
from the theory of continued fractions.
Standard references include \cite{Per29,HW,SchmLN}.

Let $\alpha := [0; a_1, a_2, \ldots]$ be a real irrational number.
Set $p_{-1}=q_0=1$ and $q_{-1}= p_0 =0$.
For $\ell \ge 1$, the $\ell$-th convergent to
$\alpha$ is the rational number $p_\ell / q_\ell := [0; a_1, a_2, \ldots , a_{\ell}]$.
Observe that
$$
q_{\ell} = a_{\ell} q_{\ell - 1} + q_{\ell-2}, \quad \ell \ge 1.   \eqno (2.1)
$$
Furthermore, the sequence $(q_{\ell})_{\ell \ge 1}$ is increasing
and $q_{\ell}$ and $q_{\ell + 1}$ are coprime for $\ell \ge 0$.

The theory of continued fraction implies that (see e.g. Theorem 164 of \cite{HW})
$$
|q_{\ell} \alpha - p_{\ell}| < q_{\ell+1}^{-1},
\quad \hbox{for $\ell \ge 1$,}   \eqno (2.2)
$$
and
$$
q_{\ell + h} \ge q_{\ell} (\sqrt{2})^{h-1},
\quad \hbox{for $h, \ell \ge 1$}. \eqno (2.3)
$$
Indeed, an easy induction on $h$ based on (2.1) proves (2.3)
for every fixed value of $\ell \ge 1$.

Likewise, an induction based on (2.1) allows us to
establish the mirror formula
$$
{q_{\ell - 1} \over q_{\ell} } = [0 ; a_\ell , a_{\ell-1}, \ldots , a_1],
\quad \ell \ge 1.    \eqno (2.4)
$$

The main tool for the proof of Theorem 1.3 is the Schmidt Subspace Theorem.

\proclaim Theorem 2.1 (W. {\bf M. Schmidt).} 
Let $m \ge 2$ be an integer.
Let $L_1, \ldots, L_m$ be linearly independent linear forms in
${\bf x} = (x_1, \ldots, x_m)$ with algebraic coefficients.
Let $\eps$ be a positive real number.
Then, the set of solutions ${\bf x} = (x_1, \ldots, x_m)$ 
in $\Z^m$ to the inequality
$$
\vert L_1 ({\bf x}) \ldots L_m ({\bf x}) \vert  \le
(\max\{|x_1|, \ldots , |x_m|\})^{-\eps}
$$
lies in finitely many proper linear subspaces of $\Q^m$.

\pro See e.g. \cite{Schm72,SchmLN}.
\cqfd

\vskip 7mm

\centerline{\bf 3. Transcendence criterion for stammering continued fractions}

\vskip 5mm

In this section we establish the part of Theorem 1.3
dealing with stammering continued fractions.
Let ${\bf a}=(a_\ell)_{\ell \ge 1}$ 
be a sequence of elements from ${\cal A}$.
We say that ${\bf a}$ 
satisfies Condition $(\spadesuit)$ if ${\bf a}$ is not
ultimately periodic and if there exist 
three sequences of finite words $(U_n)_{n \ge 1}$,   
$(V_n)_{n \ge 1}$ and $(W_n)_{n \ge 1}$ such that:

\medskip

\item{\rm (i)} For every $n \ge 1$, the word $W_n U_n V_n U_n$ is a prefix
of the word ${\bf a}$;

\smallskip

\item{\rm (ii)} The sequence
$({\vert V_n\vert} / {\vert U_n\vert})_{n \ge 1}$ is bounded from above;

\smallskip

\item{\rm (iii)} The sequence
$({\vert W_n\vert} / {\vert U_n\vert})_{n \ge 1}$ is bounded 
from above;   

\smallskip

\item{\rm (iv)} The sequence $(\vert U_n\vert)_{n \ge 1}$ is 
increasing.

\medskip

\proclaim Theorem 3.1. 
Let ${\bf a}=(a_\ell)_{\ell \ge 1}$ be a sequence of positive integers.
Let $(p_\ell/q_\ell)_{\ell \ge 1}$ denote the sequence of convergents to 
the real number
$$
\alpha:= [0; a_1, a_2, \ldots,a_\ell, \ldots].
$$ 
Assume that the sequence $(q_\ell^{1/\ell})_{\ell \ge 1}$ is bounded.
If ${\bf a}$ satisfies Condition $(\spadesuit)$, 
then $\alpha$ is transcendental.

Theorem 3.1 improves Theorem 2 from \cite{AdBu05}
and Theorem 3.1 from \cite{AdBuDa}.
Furthermore, it contains Theorem 3.2 from \cite{AdBu07e}.

Theorem 3.1 is the
exact analogue of the combinatorial transcendence criterion for expansions
to integer bases proved in \cite{AdBuLu04}.
Although its proof is very close to that of Theorem 2 of \cite{AdBu05},
we have decided to write it completely.
The new ingredient is estimate (3.4) below.

\medskip
\goodbreak
\pro
Throughout, the constants implied in $\ll$ depend only on $\alpha$.
Assume that the 
sequences $(U_n)_{n \ge 1}$,
$(V_n)_{n \ge 1}$  and $(W_n)_{n \ge 1}$ occurring in the
definition of Condition $(\spadesuit)$ are fixed.
For $n \ge 1$, set $u_n = |U_n|$, $v_n = |V_n|$ and $w_n = |W_n|$.
We assume that the real number
$\alpha:= [0; a_1, a_2, \ldots]$
is algebraic of degree at least three. 
Set $p_{-1}=q_0=1$ and $q_{-1}= p_0 =0$.

We observe
that $\alpha$ admits infinitely many good quadratic approximants
obtained by truncating its continued fraction
expansion and completing by periodicity.
Precisely, for every positive integer $n$, we define the sequence
$(b_k^{(n)})_{k \ge 1}$ by
$$
\eqalign{
b_h^{(n)} & = a_h \quad \hbox{for $1 \le h \le w_n + u_n + v_n$,} \cr
b_{w_n + h + j (u_n + v_n)}^{(n)} & = a_{w_n + h} 
\quad \hbox{for $1 \le h \le u_n + v_n$ and $j \ge 0$.} \cr}
$$
The sequence
$(b_k^{(n)})_{k \ge 1}$ is ultimately periodic, with preperiod $W_n$
and with period $U_n V_n$. Set
$$
\alpha_n= [0; b_1^{(n)}, b_2^{(n)}, \ldots , b_k^{(n)}, \ldots ]
$$
and note that, since the first $w_n + u_n + v_n + u_n$ partial
quotients of $\alpha$ and of $\alpha_n$ are the same, 
it follows from (2.2) that
$$
\Bigl| \alpha - {p_{w_n + 2 u_n + v_n} \over q_{w_n + 2 u_n + v_n}} \Bigr|
< {1 \over q_{w_n + 2 u_n + v_n}^2}
\quad \hbox{and} \quad
\Bigl| \alpha_n - {p_{w_n + 2 u_n + v_n} \over q_{w_n + 2 u_n + v_n}} \Bigr|
< {1 \over q_{w_n + 2 u_n + v_n}^2},
$$
thus,
$$
|\alpha - \alpha_n| \le  2 q_{w_n + 2 u_n + v_n}^{-2}. \eqno (3.1)
$$
Furthermore, an elementary computation (see e.g. 
\cite{Per29} on page 71) shows that $\alpha_n$ is root of the quadratic polynomial
$$
\eqalign{
P_n (X) & := (q_{w_n-1} q_{w_n+u_n + v_n} - q_{w_n} q_{w_n+u_n + v_n-1} ) X^2  \cr
&  \,  - (q_{w_n-1} p_{w_n+u_n + v_n} - q_{w_n} p_{w_n+u_n + v_n-1} 
+ p_{w_n-1} q_{w_n+u_n + v_n} - p_{w_n} q_{w_n+u_n + v_n-1}) X \cr
& \,   + (p_{w_n-1} p_{w_n+u_n + v_n} - p_{w_n} p_{w_n+u_n + v_n-1}). \cr}
$$
Since $\alpha_n$ lies in $(0, 1)$, we have $p_{\ell} \le q_{\ell}$ for $\ell \ge 1$ and
the height $H(P_n)$ of the polynomial $P_n(X)$ (the height $H(P)$
of an integer polynomial $P(X)$ is the maximum of the absolute values of its coefficients)
is at most equal to $2 q_{w_n} q_{w_n+u_n + v_n}$.
By (2.2), we have
$$
\eqalign{
& |(q_{w_n-1} q_{w_n+u_n + v_n} - q_{w_n} q_{w_n+u_n + v_n-1} ) \alpha 
- (q_{w_n-1} p_{w_n+u_n + v_n} - q_{w_n} p_{w_n+u_n + v_n-1})| \cr
& \hskip 4mm \le  \ 
q_{w_n-1} |q_{w_n+u_n + v_n} \alpha - p_{w_n + u_n + v_n}|
+ q_{w_n} |q_{w_n+u_n + v_n-1}  \alpha - p_{w_n + u_n + v_n-1}| \cr
 & \hskip 4mm \le 2  \ 
q_{w_n} \, q_{w_n + u_n + v_n}^{-1}  \cr} \eqno (3.2)
$$
and, likewise,
$$
\eqalign{
& |(q_{w_n-1} q_{w_n+u_n + v_n} - q_{w_n} q_{w_n+u_n + v_n-1} ) \alpha 
- (p_{w_n-1} q_{w_n+u_n + v_n} - p_{w_n}  q_{w_n+u_n + v_n-1})| \cr 
&  \hskip 4mmÊ\le  \ 
q_{w_n + u_n + v_n} |q_{w_n-1} \alpha - p_{w_n - 1}|
+ q_{w_n+u_n + v_n-1} |q_{w_n}  \alpha - p_{w_n}| \cr
 & \hskip 4mmÊ\le 2  \ 
 q_{w_n}^{-1} \, q_{w_n + u_n + v_n}.  \cr}  \eqno (3.3)
$$
Using (3.1), (3.2), and (3.3), we then get
$$
\eqalign{
& |P_n   (\alpha)|  = |P_n (\alpha) - P_n (\alpha_n)| \cr
& \, \, =  |(q_{w_n-1} q_{w_n+u_n + v_n} - q_{w_n} q_{w_n+u_n + v_n-1})
 (\alpha - \alpha_n) (\alpha + \alpha_n)  \cr
& - (q_{w_n-1} p_{w_n+u_n + v_n} - q_{w_n} p_{w_n+u_n + v_n-1} 
+ p_{w_n-1} q_{w_n+u_n + v_n} - p_{w_n} q_{w_n+u_n + v_n-1}) (\alpha - \alpha_n)| \cr
& \ \ = 
| (q_{w_n-1} q_{w_n+u_n + v_n} - q_{w_n} q_{w_n+u_n + v_n-1}) \alpha 
- (q_{w_n-1} p_{w_n+u_n + v_n} - q_{w_n} p_{w_n+u_n + v_n-1} ) \cr
&  + (q_{w_n-1} q_{w_n+u_n + v_n} - q_{w_n} q_{w_n+u_n + v_n-1}) 
\alpha - (p_{w_n-1} q_{w_n+u_n + v_n} - p_{w_n} q_{w_n+u_n + v_n-1}) \cr
& + (q_{w_n-1} q_{w_n+u_n + v_n} - q_{w_n} q_{w_n+u_n + v_n-1}) (\alpha_n - \alpha)| 
\cdot  |\alpha - \alpha_n| \cr
& \, \,  \ll  |\alpha - \alpha_n| \cdot \bigl( q_{w_n} \, q_{w_n + u_n + v_n}^{-1} +
q_{w_n}^{-1} \, q_{w_n + u_n + v_n} + q_{w_n} q_{w_n+u_n + v_n}  |\alpha - \alpha_n|  \bigr) \cr
& \,Ê\,  \ll  |\alpha - \alpha_n| q_{w_n}^{-1} \, q_{w_n + u_n + v_n}  \cr
& \, \,  \ll q_{w_n}^{-1} \, q_{w_n + u_n + v_n} \, q_{w_n + 2 u_n + v_n}^{-2}. \cr}  \eqno (3.4)
$$
This estimate is more precise than the upper bound
$$
|P_n (\alpha)|  \ll  H(P_n) \cdot |\alpha - \alpha_n|
\ll   q_{w_n} \, q_{w_n + u_n + v_n} \, q_{w_n + 2 u_n + v_n}^{-2}   
$$
used in \cite{AdBu05};
namely, we gain a factor $q_{w_n}^{-2}$. This improvement is crucial when $w_n$
is large.

We consider the four linearly independent linear forms:
$$
\eqalign{
L_1(X_1, X_2, X_3, X_4) = & \,  \alpha^2 X_1 - \alpha (X_2 + X_3) + X_4,  \cr
L_2(X_1, X_2, X_3, X_4) =& \, \alpha X_1 - X_2, \cr
L_3(X_1, X_2, X_3, X_4) =& \,  \alpha X_1 - X_3, \cr
L_4(X_1, X_2, X_3, X_4) =& \,  X_1. \cr}
$$
Evaluating them on the quadruple 
$$
\eqalign{
{\underline v_n} := & (q_{w_n-1} q_{w_n+u_n + v_n}  - q_{w_n} q_{w_n+u_n + v_n-1}, 
q_{w_n-1} p_{w_n+u_n + v_n} - q_{w_n} p_{w_n+u_n + v_n-1},  \cr
& p_{w_n-1} q_{w_n+u_n + v_n} - p_{w_n} q_{w_n+u_n + v_n-1}, 
p_{w_n-1} p_{w_n+u_n + v_n} - p_{w_n} p_{w_n+u_n + v_n-1}), \cr}
$$ 
it follows from (3.2), (3.3), (3.4), and (2.3) that
$$
\eqalign{
\prod_{1 \le j \le 4} \, |L_j ({\underline v_n})|
&\ll   q_{w_n + u_n + v_n}^2 \, q_{w_n + 2 u_n + v_n}^{-2} \cr
&\ll  2^{- u_n}  \cr
&\ll (q_{w_n} q_{w_n + u_n + v_n})^{- \delta u_n / (2 w_n + u_n + v_n)}, \cr}
$$
if $n$ is sufficiently large,
where we have set
$$
M =  1 + \limsup_{\ell \to + \infty} \, q_{\ell}^{1/\ell} \quad
\hbox{and} \quad \delta = {\log 2 \over \log M}.
$$
Since ${\bf a}$ satisfies Condition $(\spadesuit)$, we have
$$
\liminf_{n \to + \infty} \, {u_n \over 2 w_n + u_n + v_n}  > 0.
$$
Consequently, there exists $\eps > 0$ such that
$$
\prod_{1 \le j \le 4} \, |L_j ({\underline v_n})|
\ll ( q_{w_n} \, q_{w_n + u_n + v_n})^{-\eps}
$$
holds for any sufficiently large integer $n$.

It then follows from
Theorem 2.1 that the points ${\underline v_n}$
lie in a finite union of proper linear subspaces of $\Q^4$. 
Thus, there exist a non-zero integer quadruple $(x_1,x_2,x_3,x_4)$ and
an infinite set ${\cal N}_1$ of distinct positive integers such that
$$
\eqalign{
& x_1 (q_{w_n-1} q_{w_n+u_n + v_n} - q_{w_n} q_{w_n+u_n + v_n-1}) 
+ x_2 (q_{w_n-1} p_{w_n+u_n + v_n} - q_{w_n} p_{w_n+u_n + v_n-1})  \cr
& + x_3 (p_{w_n-1} q_{w_n+u_n + v_n} - p_{w_n} q_{w_n+u_n + v_n-1}) 
+ x_4 (p_{w_n-1} p_{w_n+u_n + v_n} - p_{w_n} p_{w_n+u_n + v_n-1}) \cr
& \, \, \, \, \,  = 0, \cr} \eqno (3.5)
$$
for any $n$ in ${\cal N}_1$.

\medskip

$\bullet$ First case: we
assume that there exist an integer $\ell$ and infinitely
many integers $n$ in ${\cal N}_1$ with $w_n = \ell$.

\medskip

By extracting an infinite subset of ${\cal N}_1$
if necessary and by considering the real number $[0; a_{\ell+1}, a_{\ell+2}, \ldots]$
instead of $\alpha$, we may without loss of generality assume
that $w_n= \ell = 0$ for any $n$ in ${\cal N}_1$. 

Then, recalling that $q_{-1}=p_0=0$ and $q_0 = p_{-1} = 1$,
we deduce from (3.5) that
$$
x_1 q_{u_n + v_n-1} + x_2  p_{ u_n + v_n-1}  - 
x_3  q_{u_n + v_n}  -  x_4   p_{u_n + v_n}  = 0,  \eqno (3.6)
$$
for any $n$ in ${\cal N}_1$.
Observe that $(x_1, x_2) \not= (0, 0)$, since, otherwise, 
by letting $n$ tend to infinity along ${\cal N}_1$ in (3.6), 
we would get that the real number
$\alpha$ is rational. Dividing (3.6) by $q_{u_n + v_n}$, we obtain
$$
x_1  {q_{u_n + v_n - 1} \over q_{u_n + v_n}} + 
x_2 { p_{u_n + v_n - 1} \over q_{u_n + v_n - 1} } \cdot {q_{u_n + v_n - 1} \over q_{u_n + v_n}}
- x_3 - x_4 {p_{u_n + v_n} \over q_{u_n + v_n}} 
= 0.  \eqno (3.7)
$$
By letting $n$ tend to infinity along ${\cal N}_1$ in (3.7),
we get that
$$
\beta := \lim_{{\cal N}_1 \ni n \to + \infty} \, {q_{u_n + v_n - 1} \over q_{u_n + v_n}}
=  {x_3 + x_4 \alpha \over x_1 + x_2 \alpha}.
$$
Furthermore, observe that, for any sufficiently large integer
$n$ in ${\cal N}_1$, we have
$$
\biggl| \beta - {q_{u_n + v_n - 1} \over q_{u_n + v_n}} \biggr| =
\biggl| {x_3 + x_4 \alpha \over x_1 + x_2 \alpha} - {x_3 + x_4 p_{u_n + v_n}/q_{u_n + v_n}
\over x_1 + x_2 p_{u_n + v_n - 1} / q_{u_n + v_n - 1}} \biggr| 
\ll {1 \over q_{u_n + v_n-1} q_{u_n + v_n}},
\eqno (3.8)
$$
by (2.2).
Since the rational number $q_{u_n + v_n-1} / q_{u_n + v_n}$ 
is under its reduced form and
$u_n + v_n$ tends to infinity when $n$ tends to infinity along ${\cal N}_1$,
we see that, for every positive real number $\eta$ and every
positive integer $N$, there exists a 
reduced rational number $a/b$ such that $b > N$ and $|\beta - a/b| \le \eta / b$.
This implies that $\beta$ is irrational.

Consider now the three linearly independent
linear forms 
$$
L'_1(Y_1, Y_2, Y_3) = \beta Y_1  - Y_2, \quad 
L'_2(Y_1, Y_2, Y_3) = \alpha Y_1 - Y_3, \quad 
L'_3(Y_1, Y_2, Y_3) = Y_2. 
$$
Evaluating them on the triple 
$(q_{u_n + v_n}, q_{u_n + v_n-1}, p_{u_n + v_n})$ with $n \in {\cal N}_1$, 
we infer from (2.2) and (3.8) that
$$
\prod_{1 \le j \le 3} \, |L'_j (q_{u_n + v_n}, q_{u_n + v_n-1}, p_{u_n + v_n})|
\ll  q_{u_n + v_n}^{-1}.
$$
It then follows from
Theorem 2.1 that the points $(q_{u_n + v_n}, q_{u_n + v_n-1}, p_{u_n + v_n})$ 
with $n \in {\cal N}_1$ lie in a finite union of proper linear subspaces of $\Q^3$. 
Thus, there exist a non-zero integer triple $(y_1, y_2, y_3)$ and
an infinite set of distinct positive integers ${\cal N}_2 \subset {\cal N}_1$ such that
$$
y_1 q_{u_n + v_n} + y_2 q_{u_n + v_n - 1} + y_3 p_{u_n + v_n}  = 0,  \eqno (3.9)
$$
for any $n$ in ${\cal N}_2$.
Dividing (3.9) by $q_{u_n + v_n}$ and letting $n$ tend to
infinity along ${\cal N}_2$, we get
$$
y_1 + y_2 \beta + y_3 \alpha = 0.  \eqno (3.10)
$$

To obtain another equation linking $\alpha$ and $\beta$, we
consider the three linearly independent
linear forms 
$$
L''_1(Z_1, Z_2, Z_3) = \beta Z_1 - Z_2, \quad 
L''_2(Z_1, Z_2, Z_3) = \alpha Z_2 - Z_3, \quad 
L''_3(Z_1, Z_2, Z_3) = Z_2. 
$$
Evaluating them on the triple 
$(q_{u_n + v_n}, q_{u_n + v_n-1}, p_{u_n + v_n-1})$ with $n$ in ${\cal N}_1$, 
we infer from (2.2) and (3.8) that
$$
\prod_{1 \le j \le 3} \, |L''_j (q_{u_n + v_n}, q_{u_n + v_n-1}, p_{u_n + v_n-1})|
\ll  q_{u_n + v_n}^{-1}.
$$
It then follows from
Theorem 2.1 that the points $(q_{u_n + v_n}, q_{u_n + v_n-1}, p_{u_n + v_n-1})$ 
with $n \in {\cal N}_1$ lie in a finite union of proper linear subspaces of $\Q^3$. 
Thus, there exist a non-zero integer triple $(z_1, z_2, z_3)$ and
an infinite set of distinct positive integers ${\cal N}_3 \subset {\cal N}_2$ such that
$$
z_1 q_{u_n + v_n} + z_2 q_{u_n + v_n - 1} + z_3 p_{u_n + v_n-1}  = 0,  \eqno (3.11)
$$
for any $n$ in ${\cal N}_3$.
Dividing (3.11) by $q_{u_n + v_n-1}$ and letting $n $ tend to
infinity along ${\cal N}_3$, we get
$$
{z_1 \over \beta} + z_2   + z_3 \alpha = 0.  \eqno (3.12)
$$
We infer from (3.10) and (3.12) that
$$
(z_3 \alpha + z_2) (y_3 \alpha + y_1) = y_2 z_1.
$$
Since $\beta$ is irrational, we get from (3.10) and (3.12) that $y_3 z_3 \not= 0$.
This shows that
$\alpha$ is an algebraic number of degree at most two, which is a contradiction
with our assumption that $\alpha$ is algebraic of degree at least three.

\medskip

$\bullet$ Second case: extracting an infinite subset ${\cal N}_4$ of ${\cal N}_1$
if necessary, we assume that $(w_n)_{n \in {\cal N}_4}$ tends to infinity.

\medskip

In particular $(p_{w_n}/q_{w_n})_{n \in {\cal N}_4}$ and 
$(p_{w_n+u_n + v_n}/q_{w_n+u_n + v_n})_{n \in {\cal N}_4}$ 
both tend to $\alpha$ as $n$ tends to infinity.

We make the following observation. Let $a$ be a letter
and $U, V, W$ be three finite words ($V$ may be empty)
such that ${\bf a}$ begins with
$W U V U$ and $a$ is the last letter of $W$ and of $U V$. 
Then, writing $W = W'a$,  $V = V'a$ if $V$ is non-empty,
and $U = U' a$ if $V$ is empty, we see
that ${\bf a}$ begins with $W' (a U) V' (a U)$ if $V$ is non-empty
and with $W' (a U') (a U')$ if $V$ is empty.
Consequently, by iterating this remark if necessary,
we can assume that
for any  $n$ in ${\cal N}_4$, the last letter of the word $U_n V_n$
differs from the last letter of the word $W_n$. Said differently, we have
$a_{w_n} \not= a_{w_n + u_n + v_n}$ for any  $n$ in ${\cal N}_4$.

Divide (3.5) by $q_{w_n} \, q_{w_n+u_n + v_n -1}$ and write 
$$
Q_n := ( q_{w_n - 1} q_{w_n+u_n + v_n} )/  (q_{w_n} q_{w_n+u_n + v_n - 1}).
$$
We then get
$$
\eqalign{
& x_1 (Q_n - 1) 
+ x_2 \biggl(Q_n { p_{w_n+u_n + v_n} \over q_{w_n+u_n + v_n}}
- {p_{w_n+u_n + v_n-1} \over q_{w_n+u_n + v_n-1} } \biggr) 
+ x_3 \biggl(Q_n { p_{w_n-1} \over q_{w_n-1}}  
- {p_{w_n} \over q_{w_n} } \biggr) \cr
& + x_4 \biggl( Q_n {p_{w_n-1} \over q_{w_n-1}}
{p_{w_n+u_n + v_n} \over q_{w_n + u_n + v_n}} - {p_{w_n} \over q_{w_n}}
{ p_{w_n+u_n + v_n-1} \over q_{w_n + u_n + v_n - 1}} \biggr) = 0, \cr} \eqno (3.13)
$$
for any $n$ in ${\cal N}_4$. To shorten the notation, for any $\ell \ge 1$,
we put $R_\ell := \alpha - p_\ell/q_\ell$ and rewrite (3.13) as
$$
\eqalign{
& \, \, \, \, \, x_1 (Q_n - 1) 
+ x_2 \bigl(Q_n (\alpha - R_{w_n+u_n + v_n})
- (\alpha - R_{w_n+u_n + v_n-1}) \bigr) \cr
& + x_3 \bigl(Q_n (\alpha - R_{w_n-1})  - (\alpha - R_{w_n}) \bigr) \cr
& + x_4 \bigl( Q_n (\alpha - R_{w_n-1})
(\alpha - R_{w_n+u_n + v_n}) - (\alpha - R_{w_n})
(\alpha -R_{w_n+u_n + v_n-1} ) \bigr) = 0. \cr} 
$$
This yields
$$
\eqalign{
& \, \, \, (Q_n-1)  \bigl( x_1 + (x_2 + x_3) \alpha + x_4 \alpha^2 \bigr) \cr
& = x_2 Q_n R_{w_n+u_n + v_n} - x_2 R_{w_n+u_n + v_n-1} + 
x_3 Q_n  R_{w_n-1} - x_3 R_{w_n} \cr
& \, \, \, - x_4 Q_n R_{w_n-1} R_{w_n+u_n + v_n}   
+ x_4 R_{w_n} R_{w_n+u_n + v_n-1} \cr 
& \, \, \, +  \alpha (x_4 Q_n R_{w_n-1} 
+ x_4 Q_n  R_{w_n+u_n + v_n} - x_4  R_{w_n} -
x_4 R_{w_n+u_n + v_n-1}). \cr} \eqno (3.14)
$$
Observe that 
$$
|R_\ell| \le q_\ell^{-1} q_{\ell+1}^{-1}, \quad \ell \ge 1, \eqno (3.15)
$$
by (2.2).

We use (3.14), (3.15) and the assumption that 
$a_{w_n} \not= a_{w_n + u_n + v_n}$ for any  $n$ in ${\cal N}_4$
to establish the following claim.

\proclaim Claim. 
We have
$$
x_1 + (x_2 + x_3) \alpha + x_4 \alpha^2 = 0.  
$$

\noindent {\it Proof of the Claim.}
If there are arbitrarily large integers $n$ in ${\cal N}_4$ such that 
$Q_n \ge 2$ or $Q_n \le 1/2$, then the claim follows from (3.14) and (3.15).

Assume that $1/2 \le Q_n \le 2$ holds for every large $n$ in ${\cal N}_4$.
We then derive from (3.14) and (3.15) that
$$
|(Q_n-1) (x_1 + (x_2 + x_3) \alpha + x_4 \alpha^2)| \ll
|R_{w_n-1}| \ll q_{w_n-1}^{-1} q_{w_n}^{-1}.
$$
If $x_1 + (x_2 + x_3) \alpha + x_4 \alpha^2 \not= 0$, then we get
$$
|Q_n - 1| \ll q_{w_n-1}^{-1} q_{w_n}^{-1}. \eqno (3.16)
$$
On the other hand, observe that, by (2.4), the rational number
$Q_n$ is the quotient of the two continued fractions
$[a_{w_n+u_n + v_n}; a_{w_n+u_n + v_n-1}, \ldots , a_1]$ and
$[a_{w_n}; a_{w_n-1}, \ldots , a_1]$. 
Since
$a_{w_n+u_n + v_n} \not= a_{w_n}$, 
we have either $a_{w_n+u_n + v_n}-a_{w_n}\geq 1$ or
$a_{w_n}-a_{w_n+u_n + v_n}\geq 1$.
In the former case, we see that
$$
Q_n \ge {a_{w_n + u_n + v_n} \over \disp a_{w_n} +
{\strut 1 \over \disp 1 + {\strut 1 \over a_{w_n - 2} + 1}}}
\ge {a_{w_n} + 1 \over \disp a_{w_n} + {a_{w_n - 2} + 1 \over a_{w_n - 2} + 2}}
\ge 1 + {1 \over (a_{w_n} + 1) (a_{w_n - 2} + 2)}.
$$
In the latter case, we have
$$
{1 \over Q_n} \ge {\disp a_{w_n} + {1 \over \disp a_{w_n - 1} + 1} \over
a_{w_n + u_n + v_n} + 1} \ge 1 + {1 \over (a_{w_n - 1} + 1) (a_{w_n + u_n + v_n} + 1)}
\ge 1 + {1 \over (a_{w_n - 1} + 1) a_{w_n}}.
$$
Consequently, in any case, we have
$$
|Q_n - 1| \gg a_{w_n}^{-1} \, \min \{  
a_{w_n - 2}^{-1},   a_{w_n - 1}^{-1} \} \gg a_{w_n}^{-1} \, q_{w_n - 1}^{-1}.
$$
Combined with (3.16), this gives 
$$
a_{w_n} \gg q_{w_n} \gg a_{w_n} q_{w_n - 1},
$$
which implies that $n$ is bounded, a contradiction.
This proves the Claim. \cqfd

Since $\alpha$ is irrational
and not quadratic, we deduce from the Claim that $x_1 = x_4 = 0$ and $x_2 = - x_3$. 
Then, $x_2$ is non-zero and, by (3.5), we have, for any $n$ in ${\cal N}_4$, 
$$
q_{w_n-1} p_{w_n+u_n + v_n} - q_{w_n} p_{w_n+u_n + v_n-1} =
p_{w_n-1} q_{w_n+u_n + v_n} - p_{w_n} q_{w_n+u_n + v_n-1}.
$$
Thus, the polynomial
$P_n(X)$ can simply be expressed as
$$
\eqalign{
& P_n (X)  := (q_{w_n-1} q_{w_n+u_n + v_n} - q_{w_n} q_{w_n+u_n + v_n-1} ) X^2  \cr
& \, \, - 2 (q_{w_n-1} p_{w_n+u_n + v_n} - q_{w_n} p_{w_n+u_n + v_n-1}) X 
+ (p_{w_n-1} p_{w_n+u_n + v_n} - p_{w_n} p_{w_n+u_n + v_n-1}). \cr}
$$

Consider now the three linearly independent linear forms 
$$
\eqalign{
L'''_1 (T_1, T_2, T_3) = & \, \alpha^2 T_1 - 2 \alpha T_2 + T_3,  \cr
L'''_2 (T_1, T_2, T_3) =  & \, \alpha  T_1 - T_2, \cr
L'''_3 (T_1, T_2, T_3) =  & \, T_1. \cr}
$$
Evaluating them on the triple 
$$
\eqalign{
{\underline v'_n} := (q_{w_n-1} q_{w_n+u_n + v_n} & - q_{w_n} q_{w_n+u_n + v_n-1}, 
q_{w_n-1} p_{w_n+u_n + v_n} - q_{w_n} p_{w_n+u_n + v_n-1},  \cr
& p_{w_n-1} p_{w_n+u_n + v_n} - p_{w_n} p_{w_n+u_n + v_n-1}), \cr}
$$ 
for $n$ in ${\cal N}_4$, it follows from (3.2) and (3.4) that
$$
\prod_{1 \le j \le 3} \, |L'''_j ({\underline v'_n})|
\ll q_{w_n}  \, q_{w_n + u_n + v_n}  \, q_{w_n +  2 u_n + v_n}^{-2} 
\ll ( q_{w_n} \, q_{w_n + u_n + v_n})^{-\eps},
$$
with the same $\eps$ as above, if $n$ is
sufficiently large.

We then deduce from
Theorem 2.1 that the points ${\underline v'_n}$, $n \in {\cal N}_4$,
lie in a finite union of proper linear subspaces of $\Q^3$. 
Thus, there exist a non-zero integer triple $(t_1, t_2, t_3)$ and
an infinite set of distinct positive integers ${\cal N}_5$
included in ${\cal N}_4$ such that
$$
\eqalign{
t_1 (q_{w_n-1} q_{w_n+u_n + v_n} & - q_{w_n} q_{w_n+u_n + v_n-1}) 
+ t_2 (q_{w_n-1} p_{w_n+u_n + v_n} - q_{w_n} p_{w_n+u_n + v_n-1})  \cr
& + t_3 (p_{w_n-1} p_{w_n+u_n + v_n} - p_{w_n} p_{w_n+u_n + v_n-1}) = 0, \cr} 
\eqno (3.17)
$$
for any $n$ in ${\cal N}_5$. 

We proceed exactly as above.
Divide (3.17) by $q_{w_n} \, q_{w_n+u_n + v_n -1}$ 
and set
$$
Q_n := ( q_{w_n - 1} q_{w_n+u_n + v_n} )/  (q_{w_n} q_{w_n+u_n + v_n - 1}).
$$
We then get 
$$
\eqalign{
t_1 (Q_n - 1) 
+ t_2 \biggl( & Q_n { p_{w_n+u_n + v_n} \over q_{w_n+u_n + v_n}}
- {p_{w_n+u_n + v_n-1} \over q_{w_n+u_n + v_n-1} } \biggr)  \cr
& + t_3 \biggl( Q_n {p_{w_n-1} \over q_{w_n-1}}
{p_{w_n+u_n + v_n} \over q_{w_n + u_n + v_n}} - {p_{w_n} \over q_{w_n}}
{ p_{w_n+u_n + v_n-1} \over q_{w_n + u_n + v_n - 1}} \biggr) = 0, \cr}  \eqno (3.18)
$$
for any $n$ in ${\cal N}_5$.
We argue as after (3.13). Since
$p_{w_n}/q_{w_n}$ and $p_{w_n+u_n + v_n}/q_{w_n+u_n + v_n}$ tend to $\alpha$
as $n$ tends to infinity along ${\cal N}_5$, we derive from (3.18) that
$$
t_1 + t_2 \alpha + t_3 \alpha^2 = 0,  
$$
a contradiction
since $\alpha$ is irrational and not quadratic. 
Consequently, $\alpha$ must be transcendental.
This concludes the proof of the theorem.   \cqfd

\vskip 7mm

\centerline{\bf 4. Proofs of Theorems 1.1 and 1.4}

\vskip 5mm

\bigskip

\noi {\it Proof of Theorem 1.1.}

Let ${\bf a} = a_1 a_2 \ldots$ be an infinite word on the alphabet $\Z_{\ge 1}$.
Assume that (1.1) does not hold. Then, there exist an integer $C \ge 2$ and 
an infinite set ${\cal N}$ of positive integers such that 
$$
p(n, {\bf a}) \le C n, \quad
\hbox{for every $n$ in ${\cal N}$}. \eqno (4.1)
$$
This implies in particular that  ${\bf a}$ is written over a finite alphabet,
thus, by (2.1), the sequence $(q_\ell^{1/\ell})_{\ell \ge 1}$ is bounded.

Let $n$ be in ${\cal N}$. By (4.1) and the
{\it Schubfachprinzip}, there exists (at least) one
block $X_n$ of length $n$ having (at least) two occurrences in the
prefix of length $(C+1)n$ of ${\bf a}$. 
Thus, there are words $W_n$, $W'_n$,
$B_n$ and $B'_n$ such that $|W_n| < |W'_n|$ and
$$
a_1 \ldots a_{(C+1)n} = W_n X_n B_n = W'_n X_n B'_n.
$$

If $|W_n X_n| \le |W'_n|$, then define $V_n$ by the equality
$W_n X_n V_n = W'_n$. Observe that
$$
a_1 \ldots a_{(C+1)n} 
= W_n X_n V_n X_n B'_n   \eqno (4.2)
$$
and
$$
{|V_n| + |W_n| \over |X_n|} \le C.   \eqno (4.3)
$$
Set $U_n :=  X_n$.

If $|W'_n| < |W_n X_n|$, then, recalling that
$|W_n| < |W'_n|$, we define $X'_n$ by
$W'_n = W_n X'_n$. Since $X_n B_n = X'_n X_n B'_n$
and $|X'_n| < |X_n|$, the word $X'_n$ is a prefix strict of $X_n$
and $X_n$ is a rational power of $X'_n$.
Thus, there are a positive integer $x_n$ and a rational number $y_n$
such that $0 \le y_n < 2$ and
$$
X'_n X_n = {X'_n}^{1 + |X_n|/|X'_n|} = 
{X'_n}^{2 x_n + y_n} = ({X'_n}^{x_n})^2 {X'_n}^{y_n}.
$$
Here and below, for a positive integer $k$, we write
$Z^k$ for the word $Z \ldots Z$ ($k$ times repeated concatenation
of the word $Z$). More generally, for any positive rational number
$r$ such that $r |Z|$ is an integer, we denote by $Z^r$ the word
$Z^{[r]} Z'$, where $Z'$ is the prefix of
$Z$ of length $(r- [r])\vert Z \vert$.  

Observe that
$$
2 x_n |X'_n| + 2 |X'_n| \ge |X'_n X_n|,
$$
thus
$$
n = |X_n| \le (2 x_n + 1) |X'_n| \le 3 x_n |X'_n|.
$$
Consequently, $W_n ({X'_n}^{x_n})^2$  is a prefix
of ${\bf a}$ such that
$$
|{X'_n}^{x_n} | \ge n/3
$$
and
$$
{|W_n| \over |{X'_n}^{x_n}|} \le {3 \over n} \cdot
\bigl( (C + 1) n - 2 |{X'_n}^{x_n}| \bigr) \le 3 C + 1. \eqno (4.4)
$$
Set $U_n := {X'_n}^{x_n}$ and let $V_n$ be the empty word.

It then follows from (4.2),
(4.3), and (4.4) that, for every $n$ in the infinite set ${\cal N}$,
$$
W_n U_n V_n U_n \quad 
\hbox{is a prefix of ${\bf a}$}
$$
with
$$
|W_n| + |V_n| \le (3 C + 1) \, |U_n|.
$$
This shows that ${\bf a}$ satisfies Condition $(\spadesuit)$.
Applying Theorem 3.1, we get that the real number
$[0; a_1, a_2, \ldots]$ is transcendental. This proves the theorem.
\cqfd

\bigskip

\noi {\it Proof of Theorem 1.4.}

If the sequence $(r_k)_{k \ge 0}$ is bounded,
then Theorem 1.4 is Corollary 3.3 of \cite{AdBu07e}.
Thus, we assume that $(r_k)_{k \ge 0}$ is unbounded and 
we consider the infinite set ${\cal K}$ composed of 
the positive integers $k$ such that $r_k > \max\{r_0, \ldots , r_{k-1}\}$.
By the assumption (1.2), there exist $\eps > 0$ and $k_0$ such that
$\lambda_{k_0} > 2$ and $\lambda_{k+1} > (1 + \eps) \lambda_k$
for $k \ge k_0$.
Let $k$ be in ${\cal K}$ with $k > k_0$.
Set
$$
W_k = a_1 a_2 \ldots a_{n_k - 1}
$$
and
$$
U_k = (a_{n_k} \ldots a_{n_k+r_k - 1})^{[\lambda_k / 2]}.
$$
Observe that ${\bf a}$ begins with $W_k U_k^2$.
Furthermore, setting
$$
n'_0 = n_0 + \sum_{h=0}^{k_0-1} \, \lambda_h r_h,
$$
we have
$$
\eqalign{
|W_k| 
& \le n'_0 +  \sum_{h=k_0}^{k-1} \, \lambda_h r_h  \cr
& \le n'_0 +  r_k \lambda_k \biggl( {1 \over 1+\eps} + \ldots 
+ {1 \over (1+\eps)^{k-k_0}} \biggr) \cr
& \le n'_0 +  r_k \lambda_k / \eps  \le 2  r_k \lambda_k / \eps  \cr}
$$
and
$$
|U_k| \ge {(\lambda_k - 1) r_k  \over  2} \ge {  \lambda_k r_k \over 4}
\ge {\eps \over 8} |W_k|,
$$
for every sufficiently large $k$ in ${\cal K}$.
Consequently, 
the word ${\bf a} = a_1 a_2 \ldots $ satisfies 
Condition $(\spadesuit)$. We conclude by applying Theorem 3.1. \cqfd

\vskip 7mm

\goodbreak

\centerline{\bf 5. Transcendence criterion for quasi-palindromic continued fractions}

\vskip 5mm

In this section, we establish the part of Theorem 1.3
dealing with quasi-palindromic continued fractions.
Let ${\bf a}=(a_\ell)_{\ell \ge 1}$ 
be a sequence of elements from ${\cal A}$.
We say that ${\bf a}$ 
satisfies Condition $(\clubsuit)$ if ${\bf a}$ is not  
ultimately periodic and if there exist 
three sequences of finite words $(U_n)_{n \ge 1}$,   
$(V_n)_{n \ge 1}$ and $(W_n)_{n \ge 1}$ such that:

\medskip

\item{\rm (i)} For every $n \ge 1$, the word $W_n U_n V_n \uU_n$ is a prefix
of the word ${\bf a}$;

\smallskip

\item{\rm (ii)} The sequence
$({\vert V_n\vert} / {\vert U_n\vert})_{n \ge 1}$ is bounded from above;

\smallskip

\item{\rm (iii)} The sequence
$({\vert W_n\vert} / {\vert U_n\vert})_{n \ge 1}$ is bounded 
from above;   

\smallskip

\item{\rm (iv)} The sequence $(\vert U_n\vert)_{n \ge 1}$ is 
increasing.

\medskip

\proclaim Theorem 5.1.
Let ${\bf a}=(a_\ell)_{\ell \ge 1}$ be a sequence of positive integers.
Let $(p_\ell/q_\ell)_{\ell \ge 1}$ denote the sequence of convergents to 
the real number
$$
\alpha:= [0; a_1, a_2, \ldots,a_\ell,\ldots].
$$ 
Assume that the sequence $(q_\ell^{1/\ell})_{\ell \ge 1}$ is bounded.
If ${\bf a}$ satisfies Condition $(\clubsuit)$,  
then $\alpha$ is transcendental.

Theorem 5.1 improves Theorem 2.4 from \cite{AdBu07f}.

\medskip

\pro
Throughout, the constants implied in $\ll$
are absolute.
We content ourselves to explain
which changes should be made to the proof of Theorem 2.4 
from \cite{AdBu07f} in order to establish Theorem 5.1. We keep
the notation of that paper.

Assume that the 
sequences $(U_n)_{n \ge 1}$, $(V_n)_{n \ge 1}$ and $(W_n)_{n \ge 1}$ are fixed. 
Set  $r_n=\vert W_n\vert$, $s_n=\vert W_n U_n\vert$ and 
$t_n=\vert W_n U_n V_n \overline{U_n}\vert$, for $n \ge 1$.
Assume that the real number 
$\alpha:= [0; a_1, a_2, \ldots]$ is algebraic of degree
at least three.

For $n \ge 1$, consider the 
rational number $P_n/Q_n$ defined by
$$
{P_n \over Q_n}:=[0;W_n U_n V_n \overline{U_n}\,\overline{W_n}]
$$
and denote by $P'_n/Q'_n$ the last convergent to $P_n/Q_n$
which is different from $P_n/Q_n$. 
It has been proved in \cite{AdBu07f} that
$$
\vert Q_n\alpha-P_n\vert <Q_nq_{t_n}^{-2}, \quad
\vert Q'_n\alpha-P'_n\vert<Q_n q_{t_n}^{-2}, \eqno (5.1)
$$
$$
\vert Q_n\alpha-Q'_n\vert<Q_n q_{s_n}^{-2},  \eqno (5.2)
$$
and
$$
Q_n \le 2 q_{r_n} q_{t_n} \le 2 q_{s_n} q_{t_n}.  \eqno (5.3)
$$
Inequality (5.2) is a consequence of the mirror formula (2.4)
which is a key ingredient for the proof of the combinatorial
transcendence criteria for quasi-palindromic continued fractions.
Since
$$
\eqalign{
\alpha(Q_n\alpha-P_n) - (Q'_n\alpha-P'_n)
& = \alpha Q_n \biggl( \alpha - {P_n \over Q_n} \biggr) -
Q'_n \biggl( \alpha - {P'_n \over Q'_n} \biggr)  \cr
& = (\alpha Q_n - Q'_n) \,  \biggl( \alpha - {P_n \over Q_n} \biggr)
+ Q'_n \biggl(  {P'_n \over Q'_n} - {P_n \over Q_n} \biggr),  \cr}
$$
it follows from (5.1), (5.2) and (5.3) that
$$
\eqalign{
|\alpha^2 Q_n - \alpha Q'_n - \alpha P_n + P'_n| & 
\ll Q_n q_{s_n}^{-2} q_{t_n}^{-2} + Q_n^{-1} \cr
& \ll Q_n^{-1}. \cr}  \eqno (5.4)
$$
Together with the 
four linearly independent linear forms with algebraic 
coefficients
$$
\eqalign{
L_1(X_1, X_2, X_3, X_4) = & \,  \alpha X_1 - X_3,  \cr
L_2(X_1, X_2, X_3, X_4) = & \, \alpha X_2 - X_4, \cr
L_3(X_1, X_2, X_3, X_4) = & \, \alpha X_1 - X_2, \cr
L_4(X_1, X_2, X_3, X_4) = & \, X_2, \cr}
$$
introduced in \cite{AdBu07f}, we consider the linear form
$$
L_5 (X_1, X_2, X_3, X_4) = 
\alpha^2 X_1 - \alpha X_2 - \alpha X_3 + X_4,
$$
and we deduce from
(5.1), (5.2), (5.3) and (5.4) that
$$
\prod_{2 \le j \le 5} \, |L_j (Q_n, Q'_n, P_n, P'_n)|
\ll Q_n^2 \, q_{t_n}^{-2} \, q_{s_n}^{-2} 
\ll q_{r_n}^2 \, q_{s_n}^{-2}.
$$
By (2.3) and (5.3), we have
$$
q_{r_n}^2 \, q_{s_n}^{-2} \ll 2^{-|U_n|} 
\ll Q_n^{-\delta(u_n + v_n - r_n)/(r_n + t_n)},
$$
if $n$ is sufficiently large,
where we have set
$$
M =  1 + \limsup_{\ell \to + \infty} \, q_{\ell}^{1/\ell} \quad
\hbox{and} \quad \delta = {\log 2 \over \log M}.
$$
Since ${\bf a}$ satisfies Condition $(\clubsuit)$,
we have
$$
\limsup_{n \to + \infty} \, {r_n \over s_n} < 1
\quad \hbox{and} \quad
\limsup_{n \to + \infty} \, {r_n + t_n \over s_n} < + \infty,
$$
thus,
$$
\liminf_{n \to + \infty} \, {u_n + v_n - r_n \over r_n + t_n} > 0.
$$
Consequently, there exists $\eps > 0$ such that
$$
\prod_{2 \le j \le 5} \, |L_j (Q_n, Q'_n, P_n, P'_n)| \ll Q_n^{-\varepsilon},   
$$
for every sufficiently large $n$.

Following the proof from \cite{AdBu07f},
we apply a first time Theorem 2.1.
It implies that the points $(Q_n, Q'_n, P_n, P'_n)$
lie in a finite union of proper linear subspaces of $\Q^4$. 
As in \cite{AdBu07f}, we deduce that there exists an
infinite set of distinct positive integers ${\cal N}$
such that $Q'_n = P_n$
for $n$ in ${\cal N}$. Thus, for $n$ in ${\cal N}$, 
we have
$$
|\alpha^2 Q_n - 2 \alpha Q'_n  + P'_n|  \ll Q_n^{-1},  \eqno (5.5)
$$
instead of (5.4).
Consider now the three linearly independent linear forms 
$$
\eqalign{
L'_1 (X_1, X_2, X_3) = & \, \alpha^2 X_1 - 2 \alpha X_2 + X_3,  \cr
L'_2 (X_1, X_2, X_3) =  & \, \alpha  X_2 - X_3, \cr
L'_3 (X_1, X_2, X_3) =  & \, X_1. \cr}
$$
Evaluating them on the triple $(Q_n, Q'_n, P'_n)$ for $n$ in ${\cal N}$,
it follows from (5.1), (5.3) and (5.5) that
$$
\prod_{1 \le j \le 3} \, |L'_j (Q_n, Q'_n, P'_n)|
\ll Q_n q_{t_n}^{-2}  \ll q_{r_n} q_{t_n}^{-1}
 \ll q_{r_n} q_{s_n}^{-1} \ll Q_n^{-\eps/2},
$$
with the same $\eps$ as above, if $n$ is
sufficiently large.

We then apply again Theorem 2.1
and we continue as in the proof of Theorem 2.4 
from \cite{AdBu07f}. We omit the details. \cqfd

\vskip 7mm

\centerline{\bf 6. Concluding remarks}

\vskip 5mm

It is likely that we are now able to get the analogues for continued
fraction expansions to all the transcendence results established recently for
expansions to an integer base and whose proofs ultimately rest
on the Schmidt Subspace Theorem.
For instance, combining the arguments of \cite{AdRa08} 
with Theorem 1.3, it is easy to prove that
if $1 \le m < M$ are integers and ${\bf a} = a_1 a_2 \ldots $ is a 
word over $\{m, M\}$ such that
$ [0; a_1, a_2, \ldots , a_{\ell}, \ldots]$ is algebraic, then 
there are arbitrarily large (finite) blocks $U$ such that $U^{7/3}$ occurs in
${\bf a}$.

Recent developments have shown that the use of
quantitative versions of the Schmidt Subspace Theorem
allows us often to strengthen or to complement results
established by means of the qualitative Schmidt Subspace Theorem; see for instance
the survey \cite{Bu11}.
In particular,
by combining ideas from \cite{AdBu10a,AdBu10c,AdBu11}
with new arguments, we have obtained in \cite{Bu12}
transcendence measures 
for transcendental real numbers whose sequence of partial quotients
${\bf a}$ is such that $n \mapsto p(n, {\bf a})/n$ is bounded.

Furthermore, proceeding as in \cite{Bu08b} and in \cite{BuEv08},
it seems to be possible to prove that if 
${\bf a} = a_1 a_2 \ldots $ is an infinite word over ${\bf Z}_{\ge 1}$ such that
$[0; a_1, a_2, \ldots , a_{\ell}, \ldots]$
is algebraic of degree at least three, then there exists $\delta > 0$ such that
$$
\limsup _{n \to + \infty} \, {p(n, {\bf a}) \over n (\log n)^{\delta}} = + \infty,   \eqno (6.1)
$$
and there exists an effectively computable positive constant $M$ such that
$$
p(n, {\bf a})  \ge \biggl( 1 + {1 \over M} \biggr) n, \quad \hbox{for $n \ge 1$}.
$$
More details will be given in a subsequent note. Observe that a statement
like (6.1) does not contain Theorem 1.1 since there exist 
infinite words ${\bf w}$ such that 
$$
\liminf_{n \to \infty} \, {p(n, {\bf w}) \over n}= 2
\quad {\rm and} \quad
\limsup_{n\to\infty}\, {p(n, {\bf w}) \over n^t}=+\infty,
\quad \hbox{for any $t>1$};
$$
see \cite{Fe96}.

\vskip 5mm

\noi {\bf Acknowledgements.} I am very thankful to the referees
for their careful reading and numerous comments, which help me to
improve the presentation of the paper.

\vskip 9mm

\bigskip

\goodbreak

\centerline{\bf References}

\vskip 5mm

\beginthebibliography{999}

\bibitem{AdBu05}
B. Adamczewski and Y. Bugeaud,
{\it On the complexity of algebraic numbers, II.
Continued fractions},
Acta Math. 195 (2005), 1--20.

\bibitem{AdBu06}
B. Adamczewski and Y. Bugeaud,
{\it Real and $p$-adic expansions involving symmetric patterns},  
Int. Math. Res. Not.  2006, Art. ID 75968, 17 pp.

\bibitem{AdBu07d}
B. Adamczewski and Y. Bugeaud,
{\it On the complexity of algebraic numbers I. 
Expansions in integer bases},
Ann. of Math. 165 (2007), 547--565.

\bibitem{AdBu07e}
B. Adamczewski and Y. Bugeaud,
{\it On the Maillet--Baker continued fractions},
J. reine angew. Math. 606 (2007), 105--121.  

\bibitem{AdBu07f}
B. Adamczewski and Y. Bugeaud,
{\it Palindromic continued fractions},
Ann. Inst. Fou\-rier (Grenoble) 57 (2007), 1557--1574.

\bibitem{AdBu10a}
B. Adamczewski et Y. Bugeaud,
{\it Mesures de transcendance et aspects quantitatifs de la 
m\'ethode de Thue--Siegel--Roth--Schmidt},
Proc. London Math. Soc. 101 (2010), 1--31.

\bibitem{AdBu10c}
B. Adamczewski and Y. Bugeaud,
{\it Transcendence measures for continued fractions involving
repetitive or symmetric patterns},
J. Europ. Math. Soc. 12 (2010), 883--914.

\bibitem{AdBu11}
B. Adamczewski et Y. Bugeaud,
{\it Nombres r\'eels de complexit\'e sous-lin\'eaire : 
mesures d'irrationalit\'e et de transcendance},
J. reine angew. Math. 658 (2011), 65--98.

\bibitem{AdBuDa}
B. Adamczewski, Y. Bugeaud, and L. Davison, 
{\it Continued fractions and transcendental numbers},
Ann. Inst. Fourier (Grenoble) 56 (2006), 2093--2113.

\bibitem{AdBuLu04}
B. Adamczewski, Y. Bugeaud et  F. Luca,
{\it Sur la complexit\'e des nombres alg\'ebriques}, 
C. R. Acad. Sci. Paris 339 (2004), 11--14.

\bibitem{AdRa08}
B. Adamczewski and N. Rampersad,
{\it On patterns occurring in binary algebraic numbers},
Proc. Amer. Math. Soc.  136  (2008),  3105--3109.

\bibitem{ADQZ}
J.-P. Allouche, J. L. Davison, M. Queff\'elec, and L. Q. Zamboni,
{\it Transcendence of Sturmian or morphic continued fractions},
J. Number Theory 91 (2001), 39--66.

\bibitem{AlSh}
J.-P. Allouche and J. Shallit, 
Automatic Sequences: Theory, Applications, Generalizations, 
Cambridge University Press, Cambridge, 2003.

\bibitem{Bak62} 
A. Baker,
{\it Continued fractions of transcendental numbers},
Mathematika 9 (1962), 1--8.


\bibitem{Bu08b}
Y. Bugeaud,
{\it An explicit lower bound for the block 
complexity of an algebraic number},
Atti Accad. Naz. Lincei Cl. Sci. Fis. 
Mat. Natur. Rend. Lincei (9) Mat. Appl. 19 (2008), 229--235.

\bibitem{Bu11}
Y. Bugeaud,
{\it Quantitative versions of the Subspace Theorem and applications},
J. Th\'eor. Nombres Bordeaux 23 (2011), 35--57.

\bibitem{Bu12}
Y. Bugeaud,
{\it Continued fractions with low complexity: 
Transcendence measures and quadratic approximation},
Compos. Math. 148 (2012), 718--750.

\bibitem{BuEv08}
Y. Bugeaud and J.-H. Evertse,
{\it On two notions of complexity of algebraic numbers},
Acta Arith. 133 (2008), 221--250.

\bibitem{Cob68}
A. Cobham,
{\it On the Hartmanis-Stearns problem for a class of tag machines}.
In: Conference Record of 1968 Ninth Annual Symposium on Switching and 
Automata Theory, Schenectady, New York (1968), 51--60.

\bibitem{Cob72}
A. Cobham,
{\it Uniform tag sequences},
Math. Systems Theory 6 (1972), 164--192.

\bibitem{Fe96}
S. Ferenczi,
{\it Rank and symbolic complexity},
Ergodic Th. Dyn. Systems 16 (1996), 663--682.

\bibitem{HW} 
G. H. Hardy and E. M. Wright, 
An introduction to the theory of numbers, 5th. edition, Clarendon Press, 1979.

\bibitem{Mai06} 
E.~Maillet, 
Introduction \`a la th\'eorie des nombres 
transcendants et des propri\'et\'es arithm\'e\-tiques des fonctions.
Gauthier-Villars, Paris, 1906.

\bibitem{MoHe38}
M. Morse and G. A. Hedlund,
{\it Symbolic dynamics},
Amer. J. Math. 60 (1938), 815--866.

\bibitem{MoHe40}
M. Morse and G. A. Hedlund,
{\it Symbolic dynamics II},
Amer. J. Math. 62 (1940), 1--42.

\bibitem{Per29} 
O.~Perron, 
Die Lehre von den Kettenbr\"uchen,
Teubner, Leibzig, 1929.

\bibitem{Que98}
M. Queff\'elec,
{\it Transcendance des fractions continues de Thue--Morse},
J. Number Theor  73 (1998), 201--211.

\bibitem{Schm72} 
W.~M.~Schmidt, 
{\it Norm form equations}, 
Ann. of Math. 96 (1972), 526--551.

\bibitem{SchmLN} 
W.~M.~Schmidt, 
Diophantine approximation,
Lecture Notes in Mathematics 785, Springer, Berlin, 1980.

\endthebibliography

\vskip1cm

\noindent Yann Bugeaud  

\noindent Universit\'e de Strasbourg

\noindent Math\'ematiques

\noindent 7, rue Ren\'e Descartes      

\noindent 67084 STRASBOURG  (FRANCE)

\vskip2mm

\noindent {\tt bugeaud@math.unistra.fr}

\bye